\documentstyle[12pt,amsmath,amssymb]{article}
\begin{document}

\newtheorem{lem}{Lemma}[section]
\newtheorem{theorem}{Theorem}[section]
\newtheorem{prop}{Proposition}[section]
\newtheorem{rem}{Remark}[section]
\newtheorem{define}{Definition}[section]
\newtheorem{cor}{Corollary}[section]
\allowdisplaybreaks

\makeatletter\@addtoreset{equation}{section}\makeatother
\def\theequation{\arabic{section}.\arabic{equation}}

\newcommand{\D}{{\cal D}}

\newcommand{\N}{{\Bbb N}}
\newcommand{\C}{{\Bbb C}}
\newcommand{\Z}{{\Bbb Z}}
\newcommand{\R}{{\Bbb R}}
\newcommand{\la}{\langle}
\newcommand{\ra}{\rangle}
\newcommand{\rom}[1]{{\rm #1}}
\newcommand{\FC}{{\cal F}C_{\mathrm b}^\infty({\cal D},\Gamma)}
\newcommand{\eps}{\varepsilon}
\newcommand{\dd}{\overset{{.}{.}}}
\newcommand{\fii}{\varphi}

\def\stackunder#1#2{\mathrel{\mathop{#2}\limits_{#1}}}
\newcommand{\FCo}{{\cal F}C_{\mathrm b}^\infty({\cal D},\dd\Gamma)}

\renewcommand{\author}[1]{\medskip{\large #1}\par\medskip}
\begin{center}{\Large \bf
  Glauber dynamics of continuous particle systems
}\end{center}

{\large Yuri Kondratiev}\\
 Fakult\"at f\"ur Mathematik, Universit\"at
Bielefeld, Postfach 10 01 31, D-33501 Bielefeld, Germany;
Institute of Mathematics, Kiev, Ukraine; BiBoS, Univ.\ Bielefeld,
Germany. e-mail:
\texttt{kondrat@mathematik.uni-bielefeld.de}\vspace{2mm}

{\large Eugene Lytvynov}\\ Institut f\"{u}r Angewandte Mathematik,
Universit\"{a}t Bonn, Wegelerstr.~6, D-53115 Bonn, Germany; BiBoS,
Univ.\ Bielefeld, Germany;  SFB 611, Univ.~Bonn, Germany. e-mail:
\texttt{lytvynov@wiener.iam.uni-bonn.de}

{\small
\begin{center}
{\bf Abstract}
\end{center}

\noindent This paper is devoted to the construction and study of
an  equilibrium Glauber-type dynamics of infinite continuous
particle systems. This dynamics is a special case of a spatial
birth and death process. On the space $\Gamma$ of all locally
finite subsets (configurations) in $\R^d$, we fix a Gibbs measure
$\mu$ corresponding to a general pair potential $\phi$ and
activity $z>0$.   We consider a Dirichlet form $ \cal E$ on
$L^2(\Gamma,\mu)$ which corresponds to the generator $H$ of the
Glauber dynamics. We  prove the existence of a Markov process $\bf
M$ on $\Gamma$ that is properly associated with $\cal E$. In the
case of a positive potential $\phi$ which satisfies
$\delta{:=}\int_{\R^d}(1-e^{-\phi(x)})\, z\, dx<1$, we also prove
that the generator $H$ has a spectral gap $\ge1-\delta$.
Furthermore, for any pure Gibbs state $\mu$, we derive a
Poincar\'e inequality.  The results about the spectral gap and the
Poincar\'e inequality are a generalization and a refinement of a
recent result from \cite{BCC}.

 \vspace{1.5mm}

\vspace{1.5mm}


\section{Introduction} This paper is devoted to the construction and study of
an equilibrium Glauber-type dynamics (GD) of infinite continuous
particle systems. This dynamics is a special case of a spatial
birth and death process on $\R^d$.  For a system of particles in a
bounded volume, such processes were introduced and studied by
C.~Preston in \cite{P}. In the latter case, the total number of
particles is finite at any moment of time.

 In the recent paper by Bertini {\it et al.},
 \cite{BCC}, the generator of the GD in
 a finite volume was studied. This generator corresponds to a special case of
 birth and death coefficients in Preston's dynamics. A positive, finite range,  pair
 potential $\phi$ and an activity $z>0$ were fixed which satisfy the condition of  the low
 activity-high temperature regime.
Then, with any finite  volume $\Lambda\subset\R^d$ and a boundary
condition $\eta$ outside $\Lambda$, one may associate a finite
volume Gibbs measure $\mu_{\Lambda,\eta}$. A non-local Dirichlet
form ${\cal
 E}_{\Lambda,\eta}$ on $L^2(\mu_{\Lambda,\eta})$ was considered which corresponds
 to the generator of the GD on $\Lambda$. It was shown  that
 the generator $H_{\Lambda,\eta}$ of  ${\cal
 E}_{\Lambda,\eta}$ has  a spectral gap which is
 uniformly positive  with respect to  all finite volumes $\Lambda$ and
 boundary conditions $\eta$.

In this paper, we discuss the GD in the infinite volume. The
problem of construction of a spatial birth and death process in
the infinite volume was initiated in paper \cite{HS}, where it was
solved in  a very special case of  nearest neighbor birth and
death processes on the real line.

So, we consider the space $\Gamma$ of all locally finite subsets
(configurations) in $\R^d$, and a grand canonical Gibbs measure
$\mu$ on $\Gamma$ which corresponds to a pair potential $\phi$ and
activity $z>0$. The measure $\mu$ is supposed to be either of the
Ruelle type or corresponding to a positive potential $\phi$
satisfying the integrability  condition.  In
Section~\ref{Section2}, we shortly recall some  facts about Gibbs
measures which we use later on.

In Section~\ref{jnhbh}, we consider the following bilinear form on
$L^2(\Gamma,\mu)$ which is defined on  a proper  set of cylinder
functions:
\begin{equation} {\cal
E}(F,G)=\int_\Gamma\sum_{x\in\gamma}(F(\gamma\setminus
x)-F(\gamma)) \,(G(\gamma\setminus
x)-G(\gamma))\,\mu(d\gamma)\label{dfhibhvi}\end{equation} (here
and below, for simplicity of notations we will just write $x$
instead of $\{x\}$ for any $x\in\R^d$). We prove that this form is
closable and its closure is a
 Dirichlet form. By using the general theory of Dirichlet forms (cf.\ \cite{MR}), we prove
 that there exists a Hunt process $\bf M$ on $\Gamma$
properly associated with $\cal E$. In particular, $\bf M$ is a
conservative Markov process on $\Gamma$ with {\it cadlag\/} paths.
By construction, $\bf M$
 is an equilibrium GD on $\Gamma$ with  the stationary measure $\mu$.
Let us mention that the birth and death coefficients were supposed
to be bounded in \cite{HS}, which is not the case for the GD,
provided the potential $\phi$ has a negative part.

In the case where the interaction between the particles is absent
(i.e., $\phi=0$ and, therefore, $\mu$ is the Poisson measure
$\pi_z$ with intensity $z$), the Markov process corresponding to
the Dirichlet form \eqref{dfhibhvi} was explicitly constructed and
studied  by D.~Surgailis \cite{S1,S2}.

In Sections~\ref{gzguz4545} and \ref{hucddhs}, we only consider
the case of a positive potential $\phi$ and study the problem  of
the spectral gap for the generator $H$ of the Dirichlet form $\cal
E$.

Let us recall that the Poisson measure $\pi_z$ possesses  the
chaos decomposition property, and hence the space
$L^2(\Gamma,\pi_z)$ is unitarily isomorphic to the symmetric Fock
space over $L^2(\R^d,z\,dx)$, see e.g.\ \cite{S1}. Under this
isomorphism, the operator $H$ goes over into the number operator
$N$ in the Fock space, see \cite[Theorem~5.1]{AKR1}. Evidently,
$N$ (and thus $H$) has  spectral gap 1. Therefore, one may expect
that, at least in the case of a ``small perturbation'' of the
Poisson measure, the operator $H$ still has a spectral gap.

One way to prove the existence of a spectral gap of a generator
$H_E$ of a Dirichlet form $ E$ is to derive a coercivity identity
for $ H_E$ on a class $\cal C$ of ``nice functions,'' and using
it, to show that, for each $F\in\cal C$, $\|H_EF\|^2\ge G (H_E
F,F)$ with $G>0$. If one additionally knows that the operator $
H_E$ is essentially selfadjoint on $\cal C$, the latter estimate
implies that $H_E$ has  a spectral gap  $\ge G$. In the case of a
probability measure defined on a Hilbert space, this approach was
developed in \cite{K}, see also \cite[Ch.~6, Sect.~4]{BK}.

So, having in mind this idea, we first prove in
Section~\ref{gzguz4545} that the operator $H$ is essentially
selfadjoint. This is technically the most difficult part of the
paper. Then, in Section~\ref{hucddhs}, we  prove a coercivity
identity for the operator $H$ on  cylinder functions, and using it
and the essential selfadjointness of $H$, we show that the set
$(0,1-\delta)$ does not belong to the spectrum of $H$, provided
that $\delta{:=}\int_{\R^d}(1-e^{-\phi(x)})\, z\, dx<1$. This
statement leads us to the Poincar\'e inequality if we are able to
show that zero is a nondegenerate eigenvalue of $H$. We prove the
latter statement  for any $\mu$ that is an extreme point in the
set of all Gibbs measures corresponding to $\phi$ and $z$. In the
low activity-high temperature regime, the latter set consists of
exactly one point, which is therefore extreme.

Thus, compared with the result of \cite{BCC}, the progress
achieved in the study of the spectral gap is as  follows:
\begin{enumerate}

\item We work in the whole space $\R^d$, instead of taking finite volumes $\Lambda$
in $\R^d$ and boundary conditions $\eta$;

\item We do not suppose that the potential $\phi$ has a finite range;

\item The essential selfadjointness of $H$ is proven;

\item For $\delta<1$, an explicit estimate for the value of the spectral gap of $H$ is
found, and a Poincar\'e inequality is proven for any pure Gibbs
state.

\end{enumerate}

In a forthcoming paper, we are going to discuss the existence
problem for  general birth and death processes on configuration
spaces and study a scaling limit of these processes.

\section{Gibbs measures on configuration spaces}\label{Section2}

The configuration space $\Gamma:=\Gamma_{\R^d}$ over $\R^d$,
$d\in\N$, is defined as the set of all subsets of $\R^d$ which are
locally finite: $$\Gamma:=\big\{\,\gamma\subset \R^d\mid
|\gamma_\Lambda|<\infty\text{ for each compact }\Lambda\subset
\R^d\,\big\},$$ where $|\cdot|$ denotes the cardinality of a set
and $\gamma_\Lambda:= \gamma\cap\Lambda$. One can identify any
$\gamma\in\Gamma$ with the positive Radon measure
$\sum_{x\in\gamma}\eps_x\in{\cal M}(\R^d)$,  where  $\eps_x$ is
the Dirac measure with mass at $x$,
$\sum_{x\in\varnothing}\varepsilon_x{:=}$zero measure, and ${\cal
M}(\R^d)$
 stands for the set of all
positive  Radon  measures on the Borel $\sigma$-algebra ${\cal
B}(\R^d)$. The space $\Gamma$ can be endowed with the relative
topology as a subset of the space ${\cal M}(\R^d)$ with the vague
topology, i.e., the weakest topology on $\Gamma$ with respect to
which  all maps $\Gamma\ni\gamma\mapsto\la
f,\gamma\ra:=\int_{\R^d} f(x)\,\gamma(dx) =\sum_{x\in\gamma}f(x)$,
$f\in{\cal D}$, are continuous. Here, $\D:=C_0(\R^d)$ is the space
of all continuous  real-valued functions on $\R^d$ with compact
support. We will denote by ${\cal B}(\Gamma)$ the Borel
$\sigma$-algebra on $\Gamma$.

Now, we proceed to consider Gibbs measures on $\Gamma$. A pair
potential  is a Borel measurable function $\phi\colon \R^d\to
{\Bbb R}\cup\{+\infty\}$ such that $\phi(-x)=\phi(x)\in\R$ for all
$x\in\R^d\setminus\{0\}$.
A grand canonical Gibbs measure $\mu $ (or just Gibbs measure for
short) corresponding to the pair potential $\phi$ and activity
$z>0$ is usually defined through the Dobrushin--Lanford--Ruelle
equation, see e.g.\ \cite{Ru69}. However, it is  convenient for us
to give an equivalent definition through the
Georgii--Nguyen--Zessin identity (\cite[Theorem~2]{NZ}, see also
\cite[Theorem~2.2.4]{Kuna}).

For $\gamma\in\Gamma$ and $x\in\R^d\setminus\gamma$, we define a
relative energy of interaction between a particle located at $x$
and the configuration $\gamma$ as follows:
 $$
E(x,\gamma){:=}\begin{cases}\sum_{y\in\gamma}\phi(x-y),&\text{if
}\sum_{y\in\gamma}|\phi(x-y)|<\infty,\\
+\infty,&\text{otherwise}.\end{cases}$$ A probability measure
$\mu$ on $(\Gamma,{\cal B}(\Gamma))$ is called a Gibbs measure if
it satisfies
\begin{equation}\int_\Gamma \mu(d\gamma)\int_{\R^d} \gamma(dx) \,
F(\gamma,x) =\int_\Gamma \mu(d\gamma)\int_{\R^d}
z\,dx\,\exp\left[-E(x,\gamma)\right] F(\gamma\cup
x,x)\label{fdrtsdrt}\end{equation} for any measurable function
$F:\Gamma\times\R^d\to[0,+\infty]$. (Notice that any fixed set
$\gamma\in\Gamma$ has zero Lebesgue measure, so that the
expression $E(x,\gamma)$ on the right hand side of
\eqref{fdrtsdrt} is a.s.\ well-defined.) Let ${\cal G}(z,\phi)$
denote the set of all Gibbs measures corresponding to $z$ and
$\phi$. In particular, if $\phi\equiv0$, then \eqref{fdrtsdrt} is
the Mecke identity, which holds if and only if $\mu$ is the
Poisson measure $\pi_z$ with intensity measure $z\, dx$.

Let us now describe some classes of  Gibbs measures which appear
in classical statistical mechanics of continuous systems.
 For every $r=(r^1,\dots,r^d)\in\Z^d$, we define
a cube $$Q_r:=\left\{\, x\in\R^d\mid r^i-\frac 12\le
x^i<r^i+\frac12 \,\right\}.$$ These cubes form a partition of
$\R^d$. For any $\gamma\in\Gamma$, we set
$\gamma_r:=\gamma_{Q_r}$, $r\in\Z^d$. For $N\in\N$ let $\Lambda_N$
be the cube with side length $2N-1$ centered at the origin in
$\R^d$, $\Lambda_N$ is then a union of $(2N-1)^d$ unit cubes of
the form $Q_r$.

 For $\Lambda\subset \R^d$, we
denote
$\Gamma_\Lambda{:=}\{\gamma\in\Gamma\mid\gamma\subset\Lambda\}$.
Now, we recall some standard conditions on $\phi$.

\begin{description}

\item[(SS)] ({\it Superstability})
 There exist $A>0$, $B\ge0$ such that, if $\gamma\in\Gamma_
{\Lambda_N}$ for some $N$, then
$$\sum_{\{x,y\}\subset\gamma}\phi(x-y)\ge\sum_{r\in\Z^d}\big(A|\gamma_r|
^2-B|\gamma_r|\big).$$

\end{description}

Notice that the superstability condition automatically implies
that the potential $\phi$ is semi-bounded from below.

\begin{description}

\item[(LR)] ({\it Lower regularity}) There exists a decreasing positive
function $a\colon\N\to{\Bbb R}_+$ such that
$$\sum_{r\in\Z^d}a(\|r\|)<\infty$$ and for any
$\Lambda',\Lambda''$ which are finite unions of cubes $Q_r$ and
disjoint, with $\gamma'\in\Gamma_{\Lambda'}$,
$\gamma''\in\Gamma_{\Lambda''}$, $$\sum_{x\in\gamma',\,
y\in\gamma'' }\phi(x-y)\ge-\sum_{r',r''\in\Z^d}a(\|r'-r''\|)
|\gamma_{r'}'|\,|\gamma_{r''}''|.$$ Here, $\|\cdot\|$ denotes the
maximum norm on $\R^d$.

\item[(I)] ({\it Integrability})
$$\int_{\R^d}|1-\exp[-\phi(x)]|\,dx<+\infty.$$

\end{description}

For results related to spectral properties of the generator of the
GD, we will need the following condition.

\begin{description}


\item[(P)] ({\it Positivity})
 $\phi(x)\ge0$ for all $x\in\R^d$.

\end{description}

A probability measure $\mu$ on $(\Gamma,{\cal B}(\Gamma))$ is
called tempered if $\mu$ is supported by
$$S_\infty{:=}\bigcup_{n=1}^\infty S_n,$$ where $$S_n:=\left\{\,
\gamma\in\Gamma\mid \forall N\in\N\ \sum_{r\in\Lambda_N\cap\Z^d}
|\gamma_r|^2\le n^2|\Lambda_N \cap\Z^d| \,\right\}.$$ By ${\cal
G}^t(z,\phi)\subset{\cal G}(z,\phi)$ we denote the set of all
tempered grand canonical Gibbs measures (Ruelle measures for
short). Due to \cite{Ru70} the set ${\cal G}^t(z,\phi)$ is
non-empty for all $z>0$ and any potential $\phi$ satisfying
conditions (SS), (LR), and (I). Furthermore, the set ${\cal
G}(z,\phi)$ is non-empty for all $z>0$ and any potential $\phi$
satisfying (P) and (I), see \cite[Proposition~2.7.15]{Kuna}.


Let us now recall  the so-called Ruelle bound (cf.\ \cite{Ru70}).

\begin{prop}\label{waessedf}  Suppose that either conditions
\rom{(I), (SS), (LR)} are satisfied and $\mu\in{\cal
G}^t(z,\phi)$\rom, $z>0$\rom, or conditions \rom{(P), (I)} are
satisfied and $\mu\in{\cal G}(z,\phi)$\rom, $z>0$\rom.  Then\rom,
for any $n\in\N$\rom, there exists a non-negative measurable
symmetric function $k_\mu^{(n)}$ on $({\Bbb R}^d)^n$
 such
that, for any measurable symmetric function
$f^{(n)}:(\R^d)^n\to[0,\infty]$,
\begin{align} &\int_\Gamma \sum_{\{x_1,\dots,x_n\}\subset\gamma}
f^{(n)}(x_1,\dots,x_n)\,\mu(d\gamma)\notag\\&\qquad =\frac1{n!}\,
\int_{(\R^d)^n} f^{(n)}(x_1,\dots,x_n)
k_\mu^{(n)}(x_1,\dots,x_n)\,dx_1\dotsm
dx_n,\label{6t565r7}\end{align} and
\begin{equation}\label{swaswea957}\forall (x_1,\dots,x_n)\in(\R^d)^n:\quad
k_\mu^{(n)}(x_1,\dots,x_n)\le \xi^n,\end{equation} where $\xi> 0$
is independent of $n$.
\end{prop}

The functions $k_\mu^{(n)}$, $n\in\N$, are called correlation
functions of the measure $\mu$, while \eqref{swaswea957} is called
the Ruelle bound.

Notice that any measure $\mu\in{\cal G}(z,\phi)$ as in
Proposition~\ref{waessedf} satisfies
\begin{equation}\label{sweqaw}\int_\Gamma \la \fii,\gamma\ra^n\,\mu(d\gamma)<\infty,\qquad
\fii\in{\cal D},\ \fii\ge0, n\in\N.\end{equation} that is,  $\mu$
has all local moments finite.

\section{The Dirichlet form $\cal E$ and associated Markov process}\label{jnhbh}

We introduce a set $\FC$ of all functions on $\Gamma$ of the form
\begin{equation}\label{hbdv}
F(\gamma)=g_F(\la\fii_1,\gamma\ra,\dots,\la\fii_N,\gamma\ra),\end{equation}
where $N\in\N$, $\fii_1,\dots,\fii_N\in\D$, and $g_F\in
C^\infty_{\mathrm b}({\Bbb R}^N)$. Here, $C^\infty_{\mathrm
b}({\Bbb R}^N)$ denotes the set of all infinitely differentiable
functions on $\R^N$ which are bounded together with all their
derivatives.  For any $\gamma\in\Gamma$, we consider
$T_\gamma{:=}L^2(\R^d,\gamma)$ as a ``tangent'' space to $\Gamma$
at the point $\gamma$, and for any $F\in\FC$ we define the
``gradient'' of $F$ at $\gamma$  as the element of $T_\gamma$
given by $D^-F(\gamma,x){:=}D^-_xF(\gamma){:=}F(\gamma\setminus
x)-F(\gamma)$, $x\in\R^d$. (Evidently, $D^-F(\gamma)$ indeed
belongs to $T_\gamma$.)

Let $\mu$ be a Gibbs measure as in Proposition~\ref{waessedf}. We
will preserve the notation $\FC$ for the set of all $\mu$-classes
of functions from $\FC$. The set $\FC$ is dense in
$L^2(\Gamma,\mu)$. We now define
\begin{align} {\cal E}(F,G){:=}&
\int_\Gamma
(D^-F(\gamma),D^-G(\gamma))_{T_\gamma}\,\mu(d\gamma)\notag
\\=&\int_\Gamma \mu(d\gamma)\int_{\R^d}\gamma(dx)\, D^-_x F(\gamma)\,
D^-_x G(\gamma), \qquad F,G\in\FC.\label{gcdtuf}\end{align} Notice
that, for any $F\in\FC$, there exists $f\in\D$ such that
$|D^-_xF(\gamma)|\le f(x)$ for all $\gamma\in\Gamma$ and
$x\in\gamma$. Hence, by \eqref{sweqaw}, the right hand side of
\eqref{gcdtuf} is well defined. By \eqref{fdrtsdrt}, we also get,
for $F,G\in\FC$,
\begin{equation}\label{iud7u}{\cal
E}(F,G)=\int_\Gamma\mu(d\gamma)\int_{\R^d} z\, dx\,
\exp\left[-E(x,\gamma)\right]D_x^+ F(\gamma)\, D_x^+
G(\gamma),\end{equation} where
$D_x^+F(\gamma){:=}F(\gamma)-F(\gamma\cup x)$.

\begin{lem}
We have\rom: ${\cal E}(F,G)=0$ for all $F,G\in\FC$ such that $F=0$
$\mu$-a\rom.e\rom.
\end{lem}

\noindent{\it Proof}. Let $F\in\FC$, $F=0$ $\mu$-a.e. Denote
$B(r){:=}\{x\in\R^d:|x|<r\}$. Then, by \eqref{fdrtsdrt},
\begin{align*} 0&=\int_\Gamma\mu(d\gamma)\int_{B(r)}\gamma(dx)\,
|F(\gamma)|\\ &=\int_\Gamma\mu(d\gamma)\int_{B(r)}z\, dx\,
\exp\left[-E(x,\gamma)\right]|F(\gamma\cup x)|.\end{align*} Hence,
$F(\gamma\cup x)=0$ for $\mu(d\gamma)\, z\,
dx\,\exp\left[-E(x,\gamma)\right]$-a.e.\
$(\gamma,x)\in\Gamma\times\R^d$. Therefore, by \eqref{iud7u}, the
lemma follows.\quad $\square$

Thus, $({\cal E},\FC)$ is a well-defined bilinear form on
$L^2(\Gamma,\mu)$.

\begin{prop} \label{ghvv}
We have\rom:
\begin{equation}\label{hvgzfc} {\cal E}(F,G)=\int_\Gamma
HF(\gamma)G(\gamma)\, \mu(d\gamma),\qquad F,G\in\FC,\end{equation}
where
\begin{equation}\label{gen}
HF(\gamma)=\int_{\R^d}´z\,
dx\,\exp\left[-E(x,\gamma)\right]D^+_xF(\gamma)-\int_{\R^d}\gamma(dx)\,D^-_xF(\gamma)
\end{equation}
and $HF\in L^2(\Gamma,\mu)$\rom. The bilinear form $({\cal
E},\FC)$ is closable on $L^2(\Gamma,\mu)$ and its closure will be
denoted by $({\cal E },D({\cal E}))$\rom. The operator $(H,\FC)$
in $L^2(\Gamma,\mu)$ has Friedrichs' extension\rom, which we
denote by $(H,D(H))$\rom.
\end{prop}

\noindent{\it Proof}. Equations \eqref{hvgzfc}, \eqref{gen} easily
follow from  \eqref{fdrtsdrt} and \eqref{gcdtuf}. Let us show that
$HF\in L^2(\Gamma,\mu)$. By \eqref{sweqaw}, the inclusion
$\int_{\R^d}\gamma(dx)\,D^-_xF(\gamma)\in L^2(\Gamma,\mu)$ is
trivial. Next, we can find a compact $\Lambda\subset \R^d$ and
$C_1>0$ such that $|D_x^+F(\gamma)|\le C_1\chi_\Lambda(x)$ for all
$\gamma\in\Gamma$ and $x\in\R^d$. Here, $\chi_\Lambda$ denotes the
indicator of $\Lambda$. Hence, by \eqref{fdrtsdrt} and
\eqref{6t565r7},
\begin{align}& \int_\Gamma\left(\int_{\R^d}´z\,
dx\,\exp\left[-E(x,\gamma)\right]D^+_xF(\gamma)\right)^2\,\mu(d\gamma)\notag\\
&\quad \le C_1^2\int_\Gamma \mu(d\gamma)\int_\Lambda z\,
dx\int_\Lambda z\, dy\,
\exp\left[-E(x,\gamma)-E(y,\gamma)-\phi(x-y)\right]\exp[\phi(x-y)]\notag\\
&\quad=2C_1^2 \int_\Gamma
\sum_{\{x,y\}\in\gamma_\Lambda}\exp[\phi(x-y)]\,\mu(d\gamma)\notag\\
&\quad= C_1^2\int_{\Lambda^2}\exp[\phi(x-y)]k_\mu^{(2)}(x,y)\,
dx\, dy<\infty,\label{zgussguz}
\end{align}
since $k_\mu^{(2)}(x,y)\le C_2\exp[-\phi(x-y)]$ for all
$x,y\in\R^d$, $C_2>0$, cf.\ \cite[Eq.~(4.29)]{AKR4}. Therefore,
$H$ is the $L^2(\Gamma,\mu)$-generator of the bilinear form ${\cal
E }$. The rest of the proposition now follows from e.g.\
\cite[Theorem~X.23]{RS}.\quad $\square$

For the notion of a ``Dirichlet form,''  appearing in the
following lemma, we refer to e.g.\ \cite[Chap.~I, Sect.~4]{MR}.

\begin{lem}\label{guzazagus}
$({\cal E},D({\cal E}))$ is a Dirichlet
form on $L^2(\Gamma,\mu)$.
\end{lem}

\noindent{\it Proof}. On $D({\cal E})$ consider the norm
$\|F\|_{D({\cal E})}{:=}\big(\|F\|^2_{L^2(\mu)}+{\cal
E}(F)\big)^{1/2}$, $F\in D({\cal E})$. Here, we denoted ${\cal E
}(F){:=}{\cal E}(F,F)$. For any $F,G\in\FC$, we define $$ {
S}(F,G)(x,\gamma){:=}D_x^-F(\gamma)\, D_x^-G(\gamma),\qquad
x\in\R^d,\ \gamma\in\Gamma.$$ Using the Cauchy--Schwarz
inequality, we conclude that $ S$ extends to a bilinear continuous
map from $(D({\cal E}),\|\cdot\|_{D({\cal E})})\times (D({\cal
E}),\|\cdot\|_{D({\cal E})})$ into
$L^1(\R^d\times\Gamma,\tilde\mu)$, where
$\tilde\mu(dx,d\gamma){:=}\linebreak\gamma(dx) \mu(d\gamma)$. Let
$F\in D({\cal E})$ and consider any sequence $F_n\in\FC$,
$n\in\N$, such that $F_n\to F$ in $(D({\cal E}),\|\cdot\|_{D({\cal
E})})$ as $n\to\infty$ . Then, $F_n(\gamma)\to F(\gamma)$ as
$n\to\infty$ for $\mu$-a.e.\ $\gamma\in\Gamma$ (if necessary, take
a subsequence of $(F_n)_{n\in\N}$ with this property).
Furthermore, for any $r>0$, we have, analogously to
\eqref{zgussguz}:
\begin{align}&
\int_{B(r)\times\Gamma }|F_n(\gamma\setminus x)-F(\gamma\setminus
x)|\, \tilde\mu(dx,d\gamma)\notag\\&\qquad =
\int_{\Gamma}\mu(d\gamma)\int_{B(r)} z\,
dx\,\exp\left[-E(x,\gamma) \right]
|F_n(\gamma)-F(\gamma)|\notag\\&\qquad\le
\left(\int_{\Gamma}|F_n(\gamma)-F(\gamma)|^2\,\mu(d\gamma)\right)^{1/2}\notag\\&\quad\qquad\times\left(\int_{\Gamma}
\left(\int_{B(r)}z\,
dx\,\exp\left[-E(x,\gamma)\right]\right)^2\,\mu(d\gamma)\right)^{1/2}\to
0\label{12345}
\end{align}
as $n\to\infty$. Therefore, $F_n(\gamma\setminus x)\to
F(\gamma\setminus x)$ for $\tilde\mu$-a.e.\
$(x,\gamma)\in\R^d\times\Gamma$. Thus, $D_x^- F_n(\gamma)\to
D_x^-F(\gamma)$ as $n\to\infty$ for $\tilde\mu$-a.e.\
$(x,\gamma)\in\R^d\times\Gamma$, which yields:
\begin{equation}\label{zguzg} { S}(F,G)(x,\gamma)=D_x^-F(\gamma)\,D^-_x
G(\gamma),\qquad\text{$\tilde\mu$-a.e.\ }(x,\gamma)\in\R^d\times
\Gamma,\ F,G\in D({\cal E}).\end{equation} Hence,
\begin{equation}\label{hh}{\cal
E}(F,G)=\int_{\R^d\times\Gamma}D^-_x F(\gamma) D_x^-G(\gamma)
\,\tilde\mu(dx,d\gamma),\qquad F,G\in D({\cal E}) .\end{equation}

Define $\R\ni x\mapsto g(x){:=}(0\vee x)\wedge 1$. Let
$(g_n)_{n\in\N}$ be a sequence of functions from $C_{\mathrm
b}^\infty (\R)$ such that, for all $x\in\R$: $0\le g_n(x)\le1$,
$0\le g_n'(x)\le 2$, $g_n(x)\to g(x)$ as $n\to\infty$.  We again
fix any $F\in D({\cal E })$ and and let $(F_n)_{n\in \N}$ be a
sequence of functions from $\FC$ such that $F_n\to F$ in $(D({\cal
E}),\|\cdot\|_{D({\cal E })})$. Consider the sequence
$(g_n(F_n))_{n\in\N}$. We evidently have: $g_n(F_n)\in\FC$ for
each $n\in\N$ and $g_n(F_n)\to g(F)$ as $n\to\infty$ in
$L^2(\Gamma,\mu)$. Next,  by the above argument,  we have:
$D_x^-g_n(F_n(\gamma))\to D_x ^-g(F(\gamma))$ as $n\to\infty$ for
$\tilde\mu$-a.e.\ $(x,\gamma)$. Furthermore, the sequence $(D^-_x
g_n(F_n(\gamma)))_{n\in\N}$ is $\tilde\mu$-uniformly
square-integrable, since so is the sequence
$(D^-_xF_n(\gamma))_{n\in\N}$. Therefore, the sequence
$(D_x^-g_n(F_n(\gamma)))_{n\in\N}$ converges to $D_x
^-g(F(\gamma))$ in $L^2(\R^d\times\Gamma,\tilde\mu)$. This yields:
$g(F)\in D({\cal E})$.

For any $x,y\in\R$, we evidently have $|g(x)-g(y)|\le|x-y|$. By
\eqref{hh}, we then finally have: ${\cal E}(g(F))\le {\cal E}(F)$,
which means that $({\cal E},D({\cal E}))$ is a Dirichlet
form.\quad $\square$

We will now need  the bigger  space $\dd\Gamma$ consisting of all
$\Z_+$-valued Radon measures on $\R^d$ (which is Polish, see e.g.\
\cite{Ka75}). Since $\Gamma\subset\dd\Gamma$ and ${\cal
B}(\dd\Gamma)\cap\Gamma={\cal B}(\Gamma)$, we can consider $\mu $
as a measure on $(\dd\Gamma,{\cal B}(\dd\Gamma))$ and
correspondingly $({\cal E}, D({\cal E}))$ as a Dirichlet form on
$L^2(\dd\Gamma,\mu)$.

For the notion of a``quasi-regular Dirichlet form,''  appearing in
the following lemma, we refer to \cite[Chap.~IV, Sect.~3]{MR}.

\begin{prop}\label{scfxgus} 
$({\cal E},D({\cal E}))$ is a
quasi-regular Dirichlet form on $L^2(\dd\Gamma,\mu)$.
\end{prop}

 \noindent {\it Proof}. By
\cite[Proposition~4.1]{MR98}, it suffices to show that there
exists a bounded, complete
 metric $\rho$ on $\dd\Gamma$ generating the vague
 topology such that, for all $\gamma_0\in\dd\Gamma$, $\rho(\cdot,\gamma_0)
\in D({\cal E})$ and $\int_{\R^d}
S(\rho(\cdot,\gamma_0))(x,\gamma)\,\gamma(dx)\le \eta(\gamma)$
$\mu$-a.e.\ for some $\eta\in L^1(\dd\Gamma,\mu)$ (independent of
$\gamma_0$). Here, $S(F){:=}S(F,F)$.

Let us recall some well-known facts about cylinder functions on
the configuration space (see e.g.\ \cite{KK} for details). Let
${\cal O}_{\mathrm c}(\R^d)$ denote the set of all open,
relatively compact sets in $\R^d$. For $\Lambda\in {\cal
O}_{\mathrm c}(\R^d)$, we  have:
$\Gamma_\Lambda=\bigsqcup_{n=0}^\infty \Gamma_\Lambda^{(n)}$,
where $\Gamma_\Lambda^{(n)}$ denotes the set of all $n$-point
subsets of $\Lambda$, $n\in\N$, and
$\Gamma_\Lambda^{(0)}=\{\varnothing\}$. For $n\in\N$, we can
naturally identify $\Gamma_\Lambda^{(n)}$ with
$\tilde\Lambda^n/S_n$, where
$\tilde\Lambda^n{:=}\{(x_1,\dots,x_n)\in\Lambda^ n: x_i\ne
x_j\text{ if }i\ne j\}$ and $S_n$ denotes the group of
permutations of $\{1,\dots,n\}$ that acts on $\tilde\Lambda^n$ by
permuting  the numbers of the coordinates. Furthermore, the trace
$\sigma$-algebra of ${\cal B}(\Gamma)$ on $\Gamma_\Lambda^{(n)}$
coincides with the $\sigma$-algebra ${\cal B}_{\mathrm sym}(\tilde
\Lambda^n )$ of all symmetric Borel subsets of $\tilde\Lambda^n$
(again under a natural isomorphism). Finally, any measurable
function $F_\Lambda$ on $\Gamma_\Lambda$ may be identified with a
measurable cylinder function $F$ on $\Gamma$ by setting
$\Gamma\ni\gamma\mapsto F(\gamma){:=}F_\Lambda(\gamma_\Lambda)$.

\begin{lem}\label{lemm1}
Let $\Lambda\in {\cal O}_{\mathrm c}(\R^d)$\rom. Any measurable
bounded function $F$ on $\Gamma_\Lambda$ such that $F\restriction
\Gamma^{(n)}_\Lambda\equiv0$ for all $n\ge N$\rom, $N\in\N$\rom,
belongs to $D({\cal E})$\rom.

\end{lem}

\noindent {\it Proof}. We take arbitrary, open, disjoint subsets
$O_1,\dots,O_n$ of $\Lambda$. Consider functions $g_1,g_2\in
C_{\mathrm b}^\infty(\R)$ such that $g_1(1)=1$, $g_2(0)=1$ and
$g_1(x)=0$ if $|x-1|>1/2$, $g_2(x)=0$ if $|x|>1/2$. Approximating
the indicator functions $\chi_{O_i}$, $i=1,\dots,n$, and
$\chi_{\Lambda\setminus (O_1\cup\dotsm\cup O_n)}$ by functions
from $\D$, we easily conclude that the statement of the lemma
holds for the function
\begin{align*}F(\gamma)&=g_1(\la\chi_{O_1},\gamma\ra)\dotsm
g_1(\la\chi_{O_n},\gamma\ra) g_2(\la \chi_{\Lambda\setminus
(O_1\cup\dotsm\cup O_n)},\gamma\ra)\\
&=\chi_{\{\,|\gamma_{O_1}|=1,\dots,|\gamma_{O_n}|=1,\,
|\gamma_{\Lambda\setminus (O_1\cup\dotsm\cup
O_n)}|=0\,\}}(\gamma).
\end{align*}
Hence, $F(\gamma)$ may be identified with the indicator function
$\chi_{S_n(O_1\times\dotsm\times O_n)}(x_1,\dots,\linebreak x_n)$
on $\tilde\Lambda^n/S_n$. Using a monotone class argument, we get
the statement for any indicator function, and then for any
measurable bounded function on $\Gamma_\Lambda^{(n)}$.\quad
$\square$

\begin{lem}\label{jjy}
Let $f:\R^d\to\R$  be measurable\rom, bounded\rom,  and with
compact support\rom. Let $\zeta\in C_{\mathrm b}^\infty(\R)$ and
$a\in\R$\rom. Then\rom, $\zeta(|\la f,\cdot\ra-a|)\in D({\cal
E})$\rom.
\end{lem}

\noindent {\it Proof}.  Let $\Lambda\in{\cal O}_c(\R^d)$ be such
that $\operatorname{supp}f\subset\Lambda$. Define
$$F_n(\gamma){:=}\zeta(|\la f,\gamma\ra-a|)\,
\chi_{\{|\gamma_\Lambda|\le n\}}(\gamma),\qquad \gamma\in\Gamma,\
n\in\N.$$ By Lemma~\ref{lemm1}, $F_n\in D({\cal E})$, $n\in\N$. We
evidently have: $F_n\to F$ in $L^2(\dd\Gamma,\mu)$. Furthermore,
using Proposition~\ref{waessedf} and the majorized convergence
theorem, we get: ${\cal E}(F_n-F_m)\to0$ as $n,m\to\infty$. From
here the statement follows. \quad $\square$

The rest of the proof of Proposition~\ref{scfxgus} is quite
analogous to the proof of \cite[Proposition~4.8]{MR98}. So, we
only outline the main changes needed.

Let $E_k{:=}B(k)\subset\R^d$, $\delta_k=1/2$,
$\gamma_0\in\dd\Gamma$ and let $g_{E_k,\delta_k}$, $\phi_k$, and
$\zeta$   be defined as in \cite{MR98}, and we additionally demand
that $\zeta'(x)\in[0,1]$ on $[0,\infty)$. Since $\phi_k
g_{E_k,\delta_k}$ is bounded and has a compact support, we have by
Lemma~\ref{jjy}: $$\zeta(|\la \phi_k g_{E_k,\delta_k},\cdot
\ra-\la\phi_k g_{E_k,\delta_k},\gamma_0\ra|)\in D({\cal E}).$$
(Notice that $\la\phi_k g_{E_k,\delta_k},\gamma_0\ra$ is a
constant.)  Furthermore, taking to notice that
$\zeta'(x)\in[0,1]$, we get from \eqref{zguzg} and the mean value
theorem:
\begin{equation} S(\zeta(|\la \phi_k g_{E_k,\delta_k},\cdot
\ra-\la\phi_k g_{E_k,\delta_k},\gamma_0\ra|)(x,\gamma)\le (\phi_k
g_{E_k,\delta_k})^2(x)\le
\chi_{B(k+1/2)}(x).\label{hahgfi}\end{equation}  Set
$$c_k{:=}\left(1+\int_{B(k+1/2)}k_\mu^{(1)}(x)\, dx
\right)^{-1/2}2^{-k/2}.$$ Using  estimate \eqref{hahgfi} and the
numbers $c_k$, we now easily obtain the statement of the
proposition absolutely analogously to the proof of
\cite[Lemma~4.11 and Proposition~4.8]{MR98}.\quad $\square$

For the notion of an ``${\cal E}$-exceptional set,'' appearing in
the next proposition, we refer to e.g.\ \cite[Chap.~III,
Sect.~2]{MR}.

\begin{prop}\label{fdj} The set $\dd\Gamma\setminus\Gamma$ is ${\cal
E}$-exceptional\rom.
\end{prop}

\noindent {\it Proof}. We modify the proof of \cite[Proposition~1
and Corollary~1]{RS98} according to our situation.

It suffices to prove the result locally, that is, to show that,
for every fixed $a\in\N$,  the set
$$N{:=}\big\{\,\gamma\in\dd\Gamma: \sup(\gamma(\{x\}): x\in
[-a,a]^d)\ge 2\,\big\}$$ is ${\cal E}$-exceptional. By
\cite[Lemma~1]{RS98}, we need to prove that there exists a
sequence $u_n\in D({\cal E})$, $n\in\N$, such that each $u_n$ is a
continuous  function on $\dd\Gamma$, $u_n\to {\bf 1} _{N}$
pointwise as $n\to\infty$, and $\sup_{n\in\N}{\cal E}
(u_n)<\infty$.

Let $f\in C_0({\Bbb R})$ be such that ${\bf 1}_{[0,1]}\le f\le
{\bf 1}_{[-1/2,3/2)}$. For any $n\in\N$ and
$i=(i_1,\dots,i_d)\in\Z^d$, define a function $f_i^{(n)} \in{\cal
D}$ by $$ f_i^{(n)}(x){:=}\prod_{k=1}^d f(n x_k- i_k),\qquad
  x\in\R^d.$$ Let also $$I_i^{(n)}(x){:=}\prod_{k=1}^d {\bf 1}
_{[-1/2,3/2)}(nx_k-i_k),\qquad x\in\R^d,$$ and note that
$f_i^{(n)}\le I_i^{(n)}$.

Let $\psi\in C_{\mathrm b}^\infty ({\Bbb R})$ be such that ${\bf
1} _{[2,\infty)}\le\psi\le{\bf 1}_{[1,\infty)}$ and $0\le\psi'\le
2\, {\bf 1}_{ (1,\infty)}$. Set ${\cal A}_n{:=}\Z^d\cap[-na,na]^d$
and define continuous functions $$\dd\Gamma\ni\gamma\mapsto
u_n(\gamma){:=}\psi \left( \sup_{i\in{\cal A}_n}\la
f_i^{(n)},\gamma\ra\right),\qquad n\in\N.$$ Evidently, $u_n\to
{\bf 1}_N$ pointwise as $n\to\infty$. Furthermore, by an
appropriate approximation of the function ${\Bbb R}^{|{\cal A}_n|}
\ni (y_1,\dots,y_{|{\cal A}_n|})\mapsto \sup_{i=1,\dots, |{\cal
A}_n|}y_i$ by $C_{\mathrm b}^\infty ({\Bbb R}^{|{\cal A}_n|})$
functions, we conclude that, for each $n\in\N$, $u_n\in{D}({\cal
E})$. We now have: $$
S(u_n)(x,\gamma)=\left(\psi\left(\sup_{i\in{\cal A }_n}\la
f_i^{(n)},\gamma-\varepsilon_x\ra
\right)-\psi\left(\sup_{i\in{\cal A }_n}\la f_i^{(n)},\gamma\ra
\right)\right)^2\quad \text{$\tilde\mu$-a.e.}$$ By the mean value
theorem, we  get,  for $\tilde\mu$-a.e.\
$(x,\gamma)\in\R^d\times\dd\Gamma $,:
\begin{equation}\label{gug}S(u_n)(x,\gamma)
=\psi'(T_n(\gamma,x))^2 \left(\sup_{i\in{\cal A }_n}\la
f_i^{(n)},\gamma-\varepsilon_x\ra-\sup_{i\in{\cal A }_n}\la
f_i^{(n)},\gamma\ra\right)^2,\end{equation} where
$T_n(\gamma,x)\in\R$ is a point between $\sup_{i\in{\cal A }_n}\la
f_i^{(n)},\gamma-\varepsilon_x\ra$ and $\sup_{i\in{\cal A }_n}\la
f_i^{(n)},\gamma\ra$. It is easy to see that the following
estimate holds, for any $\gamma\in\dd\Gamma$ and $x\in\R^d$:
\begin{align} \left|\sup_{i\in{\cal A }_n}\la
f_i^{(n)},\gamma-\varepsilon_x\ra-\sup_{i\in{\cal A }_n}\la
f_i^{(n)},\gamma\ra\right|&\le  \sup_{i\in{\cal A }_n} |\la
f_i^{(n)},\gamma-\varepsilon_x\ra-\la f_i^{(n)},\gamma\ra|\notag
\\ &=  \sup_{i\in{\cal A }_n} f_i^{(n)}(x)\notag\\
&\le \sup_{i\in{\cal A }_n} I_i^{(n)}(x)\notag \\ &\le \pmb
1_{[-a-1,a+1]^d}(x). \label{hfh}\end{align} We evidently have, for
each $\gamma\in\dd\Gamma$ and $x\in\operatorname{supp}(\gamma)$:
$$ \sup_{i\in{\cal A }_n}\la f_i^{(n)},\gamma-\varepsilon_x\ra\le
T_n(\gamma,x)\le\sup_{i\in{\cal A }_n}\la f_i^{(n)},\gamma\ra.$$
 Hence,
\begin{align}
\psi'(T_n(\gamma,x))^2&\le 4\,\pmb 1_{ \{\sup_{i\in{\cal A}_n}\la
f_i^{(n)},\cdot\ra>1\} }(\gamma)\notag\\ & \le 4\,\pmb 1_{
\{\sup_{i\in{\cal A}_n}\la I_i^{(n)},\cdot\ra\ge2\}
}(\gamma)\notag\\ &\le4\sum_{i\in{\cal A}_n}\pmb1_{\{\la
I_i^{(n)},\cdot\ra\ge2\}}(\gamma),\label{jja}
\end{align}
where we used the fact that $I_i^{(n)}$ is integer-valued. By
\eqref{gug}--\eqref{jja}, we have, for $\tilde\mu$-a.e.\
$(x,\gamma)\in\R^d\times\dd\Gamma$: $$ S(u_n)(x,\gamma)\le
4\sum_{i\in{\cal A}_n}\pmb1_{\{\la
I_i^{(n)},\cdot\ra\ge2\}}(\gamma)\, \pmb 1_{[-a-1,a+1]^d}(x).$$
Therefore, by the Cauchy--Schwarz inequality and \eqref{hh},
\begin{equation} {\cal E}(u_n)\le 4\sum_{i\in{\cal A}_n} (\mu(\{\la
I_i^{(n)},\cdot\ra\ge2\}))^{1/2}\left(\int_{\dd\Gamma}\la\gamma,\pmb1
_{[-a-1,a+1]^d}\ra\,\mu(d\gamma)\right)^{1/2}.\label{bsddbh}\end{equation}
By using \cite[Theorem~5.5]{Ru70}, we easily conclude that there
exists a constant $C_3>0 $, independent of $i$ and $n$, such that,
for all $i\in{\cal A}_n$ and $n\in\N$,
\begin{equation}\mu(\{\la I_i^{(n)},\cdot\ra\ge2\})\le C_3
\left(\int_{\R^d}I_i^{(n)}(x)\,dx\right)^2=
C_3\left(\frac2n\right)^{2d} .\label{five}
\end{equation} Since $|{\cal A}_n|=(2na+1)^d$, we get from \eqref{6t565r7},
 \eqref{bsddbh},
and \eqref{five}: $$ {\cal E}(u_n)\le 4C_3^{1/2}(2na+1)^d
\left(\frac2n\right)^d\left(\int _ {[-a-1,a+1]^d}
k_\mu^{(1)}(x)\,dx\right)^{1/2},\qquad n\in\N.$$ Therefore, there
exists a constant $C_4>0$, independent of $n$, such that ${\cal
E}(u_n)\le C_4$ for all $n\in\N$.\quad $\square$

We now have the main result of this section.

 \begin{theorem}\label{8435476} \rom{1)} Suppose that the conditions of Proposition~\rom{\ref{waessedf}}
 are satisfied\rom.  Then\rom, there
exists a Hunt process 
$${\bf M}=({\pmb{ \Omega}},{\bf F},({\bf F}_t)_{t\ge0},({\pmb
\Theta}_t)_{t\ge0}, ({\bf X}(t))_{t\ge 0},({\bf P
}_\gamma)_{\gamma\in\Gamma})$$ on $\Gamma$ \rom(see e\rom.g\rom.\
\rom{\cite[p.~92]{MR})} which is properly associated with $({\cal
E},D({\cal E}))$\rom, i\rom.e\rom{.,} for all \rom($\mu$-versions
of\/\rom) $F\in L^2(\Gamma,\mu)$ and all $t>0$ the function
\begin{equation}\label{zrd9665} \Gamma\ni\gamma\mapsto
p_tF(\gamma){:=}\int_{\pmb\Omega} F({\bf X}(t))\, d{\bf
P}_\gamma\end{equation} is an ${\cal E}$-quasi-continuous version
of $\exp(-t{H})F$\rom, where $H$ is the generator of $({\cal
E},D({\cal E}))$\rom.  $\bf M$ is up to $\mu$-equivalence unique
\rom(cf\rom.\ \rom{\cite[Chap.~IV, Sect.~6]{MR}).} In
particular\rom, $\bf M$ is $\mu$-symmetric \rom(i\rom.e\rom{.,}
$\int G\, p_tF\, d\mu=\int F \, p_t G\, d\mu$ for all
$F,G:\Gamma\to{\Bbb R}_+$\rom, ${\cal B}(\Gamma)$-measurable\rom)
and has $\mu$ as an invariant measure\rom.

\rom{2)} $\bf M$ from \rom{1)} is  up to $\mu$-equivalence
\rom(cf\rom.\ \rom{\cite[Definition~6.3]{MR}}\rom) unique between
all Hunt processes ${\bf M}'=({\pmb{ \Omega}}',{\bf F}',({\bf
F}'_t)_{t\ge0},({\pmb \Theta}'_t)_{t\ge0}, ({\bf X}'(t))_{t\ge
0},({\bf P }'_\gamma)_{\gamma\in\Gamma})$ on $\Gamma$ having $\mu$
as an invariant measure and solving the martingale problem for
$(-H, D(H))$\rom, i\rom.e\rom.\rom, for all $G\in D(H)$
$$\widetilde G({\bf X}'(t))-\widetilde G({\bf X}'(0))+\int_0^t (H
G)({\bf X}'(s))\,ds,\qquad t\ge0,$$ is an $({\bf
F}_t')$-martingale under ${\bf P}_\gamma'$ for ${\cal
E}$-q\rom.e\rom.\ $\gamma\in\Gamma$\rom. \rom(Here\rom,
$\widetilde G$ denotes a quasi-continuous version of $G$\rom,
cf\rom. \rom{\cite[Ch.~IV, Proposition~3.3]{MR}.)}
\end{theorem}

\begin{rem}\rom{In fact, the statement of Theorem~\ref{8435476}
remains true for any Gibbs measure $\mu\in{\cal G}(z,\phi)$ whose
correlation functions satisfy the Ruelle bound\rom.  }\end{rem}

\noindent {\it Proof of Theorem\/}~\ref{8435476}. The first part
of the theorem follows from
 Propositions~\ref{scfxgus}, \ref{fdj} and
\cite[Chap.~IV, Theorem~3.5 and Chap.~V, Proposition~2.15]{MR}.
The second part follows directly from (the proof of)
\cite[Theorem~3.5]{AR}.\quad $\square$

In the above theorem, $\bf M$ is canonical, i.e., $\pmb\Omega$ is
the set of all {\it cadlag} functions $\omega:[0,\infty)\to
\Gamma$ (i.e., $\omega$ is right continuous on $[0,\infty)$ and
has left limits on $(0,\infty)$), ${\bf
X}(t)(\omega){:=}\omega(t)$, $t\ge 0$, $\omega\in\pmb\Omega$,
$({\bf F}_t)_{t\ge 0}$ together with $\bf F$ is the corresponding
minimum completed admissible family (cf.\
\cite[Section~4.1]{Fu80}) and ${\pmb \Theta}_t$, $t\ge0$, are the
corresponding natural time shifts.

\section{Selfadjointness of the generator}\label{gzguz4545}

In what follows, we will always suppose that the potential $\phi$
is positive.

\begin{theorem}\label{ufhhu} Suppose that    conditions
\rom{(P), (I)} are satisfied and $\mu\in{\cal G}(z,\phi)$\rom,
$z>0$\rom. Then\rom, the operator $(H,\FC)$ is essentially
selfadjoint in $L^2(\Gamma,\mu)$\rom. In particular, Friedrichs'
extension of $(H,\FC)$ coincides with its closure\rom.
\end{theorem}

\noindent {\it Proof}. Let $(\tilde H, D(\tilde H))$ denote the
closure of $(H,\FC)$, which exists since the latter operator is a
Hermitian one. We have to show that $(\tilde H, D(\tilde H))$ is
selfadjoint. Since $\tilde H\ge0$, by the Nussbaum theorem, it is
enough to show that there exists a set ${\cal
S}\subset\bigcap_{n=1}^\infty D(\tilde H ^n)$ which is total in
$L^2(\Gamma,\mu)$ and each $F\in {\cal S}$ satisfies:
$$\sum_{n=0}^\infty \frac{\|\tilde H^n F \|_{L^2(\mu)}}{(2n)!}\,
t^n<\infty$$ for some $t>0$, see e.g.\ \cite[Theorem~X.40]{RS}.

We have the following lemma, whose proof is completely analogous
to that of Lemma~\ref{lemm1}.

\begin{lem}\label{lem1} Suppose that the conditions of
Theorem~\rom{\ref{ufhhu}} are satisfied\rom. Let a function $F$ be
as in the formulation of Lemma~\rom{\ref{lemm1}}. Then\rom, $F\in
D(\tilde H)$ and the action of $\tilde H$ on $F$ is given by the
right hand side of \eqref{gen}\rom.
\end{lem}

We denote  by ${\cal P}$ the set of all continuous polynomials on
$\Gamma$, i.e., the set of all finite sums of functions of the
form $ F(\gamma)=\prod_{i=1}^n\la \gamma,\varphi_i\ra$,  $
\varphi_i\in \D$, $i=1,\dots,n$, $n\in\N$, and constants. We
preserve the notation $\cal P$ for the set of all $\mu$-classes of
functions from $\cal P$.  Using  \eqref{swaswea957} and e.g.\
\cite{BKKL}, we see that  ${\cal P}$ is a dense subset in
$L^2(\Gamma,\mu)$. Furthermore, any function from $\cal P$ is
cylinder, and we can easily conclude from Lemma~\ref{lem1} that
its statement remains true for any function $F\in{\cal P}$.

We will now show that ${\cal P}\subset\bigcap_{n=1}^\infty
D(\tilde H^n )$. We first make some informal calculations. So, we
define:
\begin{align*} H_1F(\gamma)&{:=}\int z\, dx
\exp\left[-E(x,\gamma)\right] D^+_x F(\gamma),\\
H_2F(\gamma)&{:=}\int \gamma(dx)\, D^-_x F(\gamma),\end{align*} so
that $H=H_1-H_2$. Then,
\begin{equation}\label{fctukaf}H^n=(H_1-H_2)^n=\sum_{I\subset\{1,\dots,n\}}(-1)^{n-|I|}{\cal
H}_I^{(n)},\end{equation} where \begin{align} {\cal
H}_I^{(n)}&{:=}{\cal H}_{I,1}{\cal H}_{I,2}\dotsm {\cal
H}_{I,n},\notag \\ {\cal H}_{I,i}&{:=}\begin{cases}H_1,&i\in I,\\
H_2,&i\not\in I,\end{cases}\quad i=1,\dots,n.
\label{sdgioa}\end{align} Furthermore, by induction, we conclude:
\begin{equation}\label{whuwiuiw} {\cal H}_I^{(n)}=\sum_{J\subset \big\{1,\dots,\,
\max\{i:\, i\in I\}-1\big\}} {\cal H}_{I,J}^{(n)},\end{equation}
where
\begin{align}
&{\cal H}_{I,J}^{(n)}F(\gamma)=\bigg(\int m_{I,n}(dx_n)\,
U_{I,J,n,x_n}\int m_{I,n-1}(dx_{n-1})U_{I,J,n-1,x_{n-1}}\notag
\\&\qquad  \dotsm \int m_{I,1}(dx_1)\, U_{I,J,1,x_1}
 \prod_{i\in I}\exp\left[-\sum_{u\in\eta\setminus\{x_s:\,
s\in I^c,\, s>i\}}  \phi(x_i-u)\right] \notag \\ &\qquad \times
\prod_{j=1}^n
G_{I,J,j}(x_1,\dots,x_n)F(\gamma)\bigg){\bigg|}_{\eta=\gamma},\notag\\
& I^c{:=}\{1,\dots,n\}\setminus I,\notag\\ &
m_{I,i}(dx_i){:=}\begin{cases}z\, dx_i,&i\in I,\\ \gamma(dx_i),&
i\in I^c,\end{cases}\notag\\ &
U_{I,J,i,x_i}{:=}\begin{cases}D^+_{x_i},& i\in I,\ i\in J^c,\\
D^-_{x_i},& i\in I^c,\ i\in J^c,\\ \operatorname{id},& i\in
J,\end{cases}\notag \\ &
G_{I,J,j}^{(n)}(x_1,\dots,x_n){:=}\begin{cases}\exp\left[-\sum_{r\in
I,\, r<j }\phi(x_j-x_r)\right],& j\in I,\ j\in J^c,\\
1-\exp\left[-\sum_{r\in I,\, r<j }\phi(x_j-x_r)\right],& j\in J,\\
1,& j\in I^c,\ j\in J^c,\end{cases}\notag\\ &i,j=1,\dots,n.
\label{sigidz}
\end{align}
For example, let $n=9$, $I=\{3,4,5,8\}$, $J=\{6,7\}$. Then,
\begin{align*}
{\cal H}_{I,J}^{(n)}F(\gamma)&=\bigg( \int\gamma(dx_9)\,D^-_{x_9}
\int z\, dx_8\, D^+_{x_8}\int\gamma(dx_7)\int\gamma(dx_6)\int z\,
dx_5\, D^+_{x_5}\\&\quad \times \int z\, dx_4\, D^+_{x_4} \int z\,
dx_3\, D_{x_3}^+ \int\gamma(dx_2)\,
D^-_{x_2}\int\gamma(dx_1)\,D^-_{x_1}\\&\quad \times
\exp\left[-\sum_{u\in\eta\setminus\{x_9\}}\phi(x_8-u)\right]\\
&\quad \times \exp\left[-
\sum_{u\in\eta\setminus\{x_6,x_7,x_9\}}\big(
\phi(x_5-u)+\phi(x_4-u)+\phi(x_3-u)\big)\right]\\ &\quad \times
\exp\big[-\phi(x_8-x_5)-\phi(x_8-x_4)-\phi(x_8-x_3)\big]\\&\quad\times
\big(1-\exp\big[-\phi(x_7-x_5)-\phi(x_7-x_4)-\phi(x_7-x_3)\big]\big)
\\&\quad\times
\big(1-\exp\big[-\phi(x_6-x_5)-\phi(x_6-x_4)-\phi(x_6-x_3)\big]\big)\\
&\quad\times
\exp\big[-\phi(x_5-x_4)-\phi(x_5-x_3)-\phi(x_4-x_3)\big]
F(\gamma)\bigg){\bigg|}_{\eta=\gamma}.
\end{align*}

\begin{lem}\label{lem2} Suppose that the conditions of
Theorem~\rom{\ref{ufhhu}} are satisfied\rom.  We then have ${\cal
P}\subset\bigcap_{n=1}^\infty D(\tilde H^n)$\rom, and for any
$F\in{\cal P}$\rom, $\tilde H^n F$ is given by formulas
\eqref{fctukaf}--\eqref{sigidz} \rom(in which $H$ is replaced by
$\tilde H$\rom)\rom.
\end{lem}

\noindent {\it Proof}. This statement follows from
\eqref{swaswea957}, Lemma~\ref{lem1} and formulas
\eqref{fctukaf}--\eqref{sigidz}. Indeed, replace in formulas
\eqref{fctukaf}--\eqref{sigidz} $\phi(x)$ by the function
$\phi_n(x){:=} \phi(x) \chi_{B(n)}(x)$, $n\in\N$, and take $F\in
\cal P$. Then, the obtained functions become cylindrical.
Approximate these by functions as in Lemma~\ref{lem1} and let
$n\to\infty$. The rest then easily follows.\quad $\square$

\begin{lem}\label{lem3} Suppose that the conditions of
Theorem~\rom{\ref{ufhhu}} are satisfied\rom. Then\rom, for any $
F(\gamma)=\prod_{i=1}^l\la \gamma,\varphi_i\ra$,  $ \varphi_i\in
\D$, $i=1,\dots,l$, $l\in\N$\rom, there exists $t>0$ such that
$\sum_{n=0}^\infty \|\tilde H^n F
\|_{L^2(\mu)}t^n/(2n)!<\infty$\rom.
\end{lem}

\noindent {\it Proof}. We first derive some estimates.

Let $f_1,\dots,f_k$ be bounded integrable functions on $\R^d$.
Consider $G(\gamma){:=}\linebreak\prod_{i=1}^k\la \gamma, f_i\ra$.
From \eqref{6t565r7}, we conclude that
\begin{align*}&\int_\Gamma G(\gamma)\,\mu(d\gamma)= \sum_{i=1}^k
\sum_{\stackunder{\text{$A_j$'s disjoint},\, A_1\cup\dotsm\cup
A_i=\{1,\dots,k\} }{(A_1,\dots,A_i):\, \varnothing \ne
A_j\subset\{1,\dots,k\},\, j=1,\dots,i}}\\&\qquad\times
\int_{(\R^d)^i} g^{(k)}(\underbrace{x_1,\dots,x_1}_{\text{$|A_1|$
times}},\underbrace{x_2,\dots,x_2}_{\text{$|A_2|$ times}},\dots,
\underbrace{x_i,\dots,x_i}_{\text{$|A_i|$ times}})
k_\mu^{(i)}(x_1,\dots,x_i)\, dx_1\dotsm dx_i,
\end{align*}
where $g^{(k)}(x_1,\dots,x_k)=(1/k!)\sum_{\sigma\in
S_k}f_1(x_{\sigma(1)})\dotsm f_k(x_{\sigma(k)})$. By induction, we
prove $$
 \sum_{i=1}^k
\sum_{\stackunder{\text{$A_j$'s disjoint},\, A_1\cup\dotsm\cup
A_i=\{1,\dots,k\} }{(A_1,\dots,A_i):\, \varnothing \ne
A_j\subset\{1,\dots,k\},\, j=1,\dots,i}} 1\le 2^{k-1}k!,\qquad
k\in\N.$$ Therefore, by \eqref{swaswea957},
\begin{equation}\label{a}\int_\Gamma |G(\gamma)|\,\mu(d\gamma)\le
2^ {k-1}\max\{1,\xi\}^k\, k!\,
\prod_{i=1}^k\max\{\|f_i\|_{L^1},\|f_i\|_{L^\infty}\}.\end{equation}
Note that, for $\mu$-a.e.\ $\gamma\in\Gamma$,
\begin{align}D_x^+ G(\gamma)&=-\sum_{i=0}^{k-1} \binom i k \int \gamma
(dy_1)\dotsm
\gamma(dy_{i})\,g^{(k)}(y_1,\dots,y_{i},\underbrace{x,\dots,x}_{\text{$(k-i)$
times }}),\notag\\ D^-_x G(\gamma)&=\sum_{i=0}^{k-1} \binom i k
(-1)^{k-i}\int \gamma (dy_1)\dotsm
\gamma(dy_{i})\,g^{(k)}(y_1,\dots,y_{i},\underbrace{x,\dots,x}_{\text{$(k-i)$
times }}).\label{c}
\end{align}
We next easily get the following identity:
\begin{equation}\label{d}
\sum_{k_1=0}^n\binom
{k_1}n\sum_{k_2=0}^{k_1}\binom{k_2}{k_1}\sum_{k_3=0}^{k_2}\binom{k_3}{k_2}\dotsm
\sum_{k_m=0}^{k_{m-1}}\binom{k_m}{k_{m-1}}=m^n,\qquad m,n\in\N.
\end{equation}
Recall also the standard estimate
\begin{equation}\label{e}(a_1+\dots+a_n)^2\le
n(a_1^2+\dots+a_n^2),\qquad a_1,\dots,a_n\in\R,\
n\in\N.\end{equation} Finally, using condition (P), we get, for
any $x,y_1,\dots,y_k\in\R^d$ and $k\in\N$,
\begin{equation}\label{f}
1-\exp\left[-\sum_{i=1}^k\phi(x-y_i)\right]\le\sum_{i=1}^k(1-\exp[-\phi(x-y_i)]).
\end{equation}

Let now $F(\gamma)$ be as in the formulation of the lemma. Fix
$\Lambda\in {\cal O}_{\mathrm c}(\R^d)$ and $C_5>0$ such that
\begin{equation}\label{g}
|f_1(x_1)\dotsm f_l(x_l)|\le
C_5\chi_{\Lambda^l}(x_1,\dots,x_l),\qquad x_1,\dots,x_l\in\R^d.
\end{equation}
For any sets $I,J\subset\{1,\dots,n\}$ as in \eqref{whuwiuiw}, we
define $$ n_1{:=}|I\cap J|,\quad n_2{:=}|I^c\cap J|,\quad
n_3{:=}|I\cap J^c|,\quad n_4{:=}|I^c\cap J^c|.$$ Notice that
$n_1+n_2+n_3+n_4=n$. Estimating all the multipliers of the form
$e^{-\phi(\cdot)}$ by 1, and using \eqref{a}--\eqref{g}, we get
from \eqref{sigidz} \begin{align}& \|{\cal
H}_{I,J}^{(n)}F\|_{L^2(\mu)}^2\le \max\{1,z\}^{2n}C_5^2
\max\left\{1,\int_\Lambda
dx\right\}^{2(n+l)}\max\{1,\xi\}^{2(n+l)}\notag\\ &\qquad\times
\max\left\{1,\int_{\R^d}(1-\exp[-\phi(x)])\right\}^{2n}\notag\\
&\qquad\times 2^{2(l+n_2)-1}(2(l+n_2))!\,(n_1+n_3)!\,
(n_1+n_3)^{n_2}(n_3+n_4)^{l+n_2},\label{vydvsj}
\end{align}
where the factor $2^{2(l+n_2)-1}(2(l+n_2))!$ is connected with
estimate \eqref{a} and the fact that we get monomials of order
$\le l+n_2$, the factor $(n_1+n_3)!$ is connected with the
application of \eqref{f} to the terms connected with the set $I$,
the factor $(n_1+n_3)^{n_2}$ is connected with the application of
\eqref{f} to the terms connected with $I^c\cap J$, and the factor
$(n_3+n_4)^{l+n_2}$ is connected with the application of \eqref c,
\eqref d to a monomial of order $\le(l+n_2)$. Using the estimate
$(2k)!\le 4^k(k!)^2$, $k\in\N$, we conclude from \eqref{vydvsj}
that there exists a constant $C_6>0$, independent of $n,I,J$ and
thus depending only on $F$, such that $$ \|{\cal
H}_{I,J}^{(n)}F\|_{L^2(\mu)}^2\le C_6^n (n!)^2 n^{2n}.$$ Hence, by
\eqref{fctukaf} \begin{equation}\label{dcgh} \|\tilde H^n
f\|_{L^2(\mu)}\le (2C_6)^{n/2}n!\, n^n.\end{equation} Estimate
\eqref{dcgh}, together with Stirling's formula, easily implies the
statement of the lemma, and hence the statement of the
theorem.\quad$\square$

\section{Spectral gap of the generator}\label{hucddhs}

We first prove a coercivity identity for the gradient $D^-$. We
note that, for any $\gamma\in\Gamma$ and $F\in\FC$, $(D^-)^2
F(\gamma)$ is the element of the Hilbert space $T_\gamma^{\otimes
2}=L^2((\R^d)^2,\gamma^{\otimes 2})$ given by
$(D^-)^2F(\gamma,x,y)=D^-_xD^-_yF(\gamma)$, $x,y\in\R^d$.
Furthermore, for any $x,y\in\gamma$: \begin{equation}\label{hsg}
D^-_xD^-_y F(\gamma)=\begin{cases}
F(\gamma\setminus\{x,y\})-F(\gamma\setminus x)-F(\gamma\setminus
y)+F(\gamma),&x\ne y,\\ F(\gamma)-F(\gamma\setminus
x)=-D^-_xF(\gamma),& x=y.\end{cases}\end{equation} Through the
natural identification of the elements of $T_\gamma^{\otimes 2}$
with linear continuous operators in $T_\gamma$, we get:
\begin{equation}\label{uisdz}\operatorname{Tr}
(D^-)^2F(\gamma)((D^-)^2F(\gamma))^*=\sum_{x,y\in\gamma}(D^-_xD^-_yF(\gamma))^2.\end{equation}
Here, $\operatorname{Tr}$ denotes the trace of an operator and
$((D^-)^2F(\gamma))^*$ is the adjoint operator of
$(D^-)^2F(\gamma)$.

\begin{lem}[Coercivity identity]\label{zsfkllk} Suppose that the conditions of
Theorem~\rom{\ref{ufhhu}} are satisfied\rom. Then\rom, for any
$F\in\FC$\rom, we have
\begin{align*}&\int_\Gamma
(HF(\gamma))^2\,\mu(d\gamma)=\int_\Gamma\bigg[
\operatorname{Tr}(D^-)^2F(\gamma)((D^-)^2F(\gamma))^*
+\sum_{x,y\in\gamma,\, x\ne
y}(\exp[\phi(x-y)]-1)\\&\qquad\times(F(\gamma\setminus\{x,y\})-F(\gamma\setminus
x)) (F(\gamma\setminus\{x,y\})-F(\gamma\setminus
y))\bigg]\,\mu(d\gamma).\end{align*}
\end{lem}

\noindent{\it Proof}. Analogously to \eqref{zgussguz}, we get from
\eqref{fdrtsdrt}:
\begin{align}&\int_\Gamma \mu(d\gamma)\,\left(\int_{\R^d}z \, dx\, \exp\left[-
E(x,\gamma)\right]D_x^+F(\gamma)\right)^2\notag\\&\qquad =
\int_\Gamma \mu(d\gamma)\sum_{x,y\in\gamma,\,x\ne
y}\exp[\phi(x-y)](F(\gamma\setminus\{x,y\})-F(\gamma\setminus
x))\notag\\&\qquad\quad\times
(F(\gamma\setminus\{x,y\})-F(\gamma\setminus
y)).\label{gaft}\end{align} Next,
\begin{align}&
-2\int_\Gamma\mu(d\gamma)\int_{\R^d}z\, dx\,
\exp\left[-E(x,\gamma)\right]D_x^+
F(\gamma)\int_{\R^d}\gamma(dy)\,D^-_y F(\gamma)\notag\\ &\qquad=
-2\int_\Gamma\mu(d\gamma)\int_{\R^d}\gamma(dx)\,(F(\gamma\setminus
x)-F(\gamma))\notag\\&\qquad\quad\times\int_{\R^d}
(\gamma\setminus
x)(dy)\,(F(\gamma\setminus\{x,y\})-F(\gamma\setminus
x))\notag\\&\qquad=-
\int_\Gamma\mu(d\gamma)\sum_{x,y\in\gamma,\,x\ne
y}\big[(F(\gamma\setminus
x)-F(\gamma))(F(\gamma\setminus\{x,y\})-F(\gamma\setminus
x))\notag\\ &\qquad\quad- (F(\gamma\setminus
y)-F(\gamma))(F(\gamma\setminus\{x,y\})-F(\gamma\setminus
y))\big].\label{hdh}
\end{align}
Finally,
\begin{align}&\int_\Gamma\mu(d\gamma)\left(\int_{\R^d}\gamma(dx)\,
D^-_xF(\gamma)\right)^2\notag\\&\qquad=\int_\Gamma
\mu(d\gamma)\sum_{x\in\gamma}(F(\gamma\setminus
x)-F(\gamma))^2\notag\\&\qquad\quad +\int_\Gamma
\mu(d\gamma)\sum_{x,y\in\gamma,\,x\ne y}(F(\gamma\setminus
x)-F(\gamma))(F(\gamma\setminus
y)-F(\gamma))\label{ghjs}\end{align} By \eqref{gen} and
\eqref{hsg}--\eqref{ghjs}, the lemma follows.\quad $\square$

\begin{theorem}\label{ewsrth080} Suppose that \rom{(P)} holds\rom, $z>0$\rom,
and
\begin{equation} \delta{:=}\int_{\R^d}(1-\exp[-\phi(x)])\, z\,
dx<1.\label{guzggz}
\end{equation} Let $\mu\in{\cal G}(z,\phi)$\rom. Then\rom, the set $(0,1-\delta)$ does not belong to the
spectrum of $H$\rom.
\end{theorem}

\noindent {\it Proof}. We fix any $F\in\FC$. By \eqref{hsg} and
\eqref{uisdz}, we have:
\begin{equation}\label{jvh}\operatorname{Tr}
(D^-)^2F(\gamma)((D^-)^2F(\gamma))^*
\ge\sum_{x\in\gamma}(D^-_xD^-_xF(\gamma))^2=\sum_{x\in\gamma}(D^-_xF(\gamma))^2,\qquad
\gamma\in\Gamma.\end{equation} Using (P), \eqref{fdrtsdrt},
\eqref{guzggz}, and the Cauchy--Schwarz inequality, we next have
\begin{align} &\bigg|\int_\Gamma \sum_{x,y\in\gamma,\, x\ne
y}(\exp[\phi(x-y)]-1)(F(\gamma\setminus\{x,y\})-F(\gamma\setminus
x))\notag\\ &\qquad\quad\times
(F(\gamma\setminus\{x,y\})-F(\gamma\setminus
y))\,\mu(d\gamma)\bigg|\notag\\ &\qquad \le \int_\Gamma
\sum_{x,y\in\gamma,\, x\ne
y}(\exp[\phi(x-y)]-1)(F(\gamma\setminus\{x,y\})-F(\gamma\setminus
y))^2 \,\mu(d\gamma)\notag\\ &\qquad=
\int_\Gamma\mu(d\gamma)\int_{\R^d}\gamma(dy)\int_{\R^d}(\gamma\setminus
y)(dx)\, (\exp[\phi(x-y)]-1)\notag\\ &\qquad\quad\times
(F(\gamma\setminus\{x,y\})-F(\gamma\setminus y))^2\notag\\
&\qquad=\int_\Gamma\mu(d\gamma)\int_{\R^d}z\,
dy\,\exp\left[-E(y,\gamma)\right]\int_{\R^d}\gamma(dx)\,(\exp[\phi(x-y)]-1)\notag\\
&\qquad\quad\times (F(\gamma\setminus x)-F(\gamma))^2\notag\\
 &\qquad=
\int_\Gamma\mu(d\gamma)\int_{\R^d}z\, dy\int_{\R^d}\gamma(dx)\,
\exp\left[-E(y,\gamma\setminus  x  )\right]\notag\\
&\qquad\quad\times(1-\exp[-\phi(x-y)])(F(\gamma\setminus
x)-F(\gamma))^2\notag\\ &\qquad\le
\int_\Gamma\mu(d\gamma)\int_{\R^d}z\, dy\int_{\R^d}\gamma(dx)\,
(1-\exp[-\phi(x-y)])(F(\gamma\setminus
x)-F(\gamma))^2\notag\\&\qquad= \int_{\R^d}(1-\exp[-\phi(y)])\,
z\,dy\times
\int_\Gamma\mu(d\gamma)\int_{\R^d}\gamma(dx)\,(D^-_xF(\gamma))^2\notag\\&\qquad=
\delta\, (HF,F)_{L^2(\mu)}.\label{huui}
\end{align}
Using Lemma~\ref{zsfkllk}, \eqref{jvh}, and \eqref{huui}, we get,
for each $F\in\FC$:
\begin{equation}
(HF,HF)_{L^2(\mu)}\ge(1-\delta)(HF,F)_{L^2(\mu)}.\label{shjn}\end{equation}
From Theorem~\ref{ufhhu}, we then conclude that \eqref{shjn} holds
true for each $F\in D(H)$. Therefore, denoting by
$(E_\lambda)_{\lambda\ge0}$ the resolution of the identity of the
operator $H$, we have:
$$\int_{[0,\infty)}\lambda(\lambda-(1-\delta))\, d(E_\lambda
F,F)_{L^2(\mu)}\ge0,\qquad F\in D(H). $$ From here the statement
of the theorem trivially follows.\quad $\square$

\begin{cor}[Poincar\'e inequality]\label{jndfjdj4}
Suppose  \rom{(P)} and \rom{\eqref{guzggz}} hold and suppose $\mu$
is an extreme point of the convex set ${\cal G}(z,\phi)$\rom.
Then\rom,
\begin{equation}
{\cal E}(F,F)\ge(1-\delta)\int_\Gamma(F(\gamma)-\la
F\ra_\mu)^2\,\mu(d\gamma),\qquad F\in D({\cal E}),
\label{chdhbbh}\end{equation} where $\la F\ra_\mu{:=}\int_\Gamma
F(\gamma)\,\mu(d\gamma)$.
\end{cor}

\begin{rem}
\rom{The Poincar\'e inequality \eqref{chdhbbh} means that,  in
addition to the fact that the set $(0,1-\delta)$ does not belong
to the spectrum of $H$, we also have that the kernel of $H$
consists only of the constants. }
\end{rem}

\noindent {\it Proof of Corrolary~\rom{\ref{jndfjdj4}}}. Since
$\mu$ is extreme in ${\cal G}(z,\phi)$, analogously to proof of
the part (i)$\Rightarrow$(ii) of \cite[Theorem~6.2]{AKR4}, we
conclude:
\begin{align}& \{ \nu\in {\cal G}(z,\phi) \mid \nu=\rho\cdot \mu
\text{\rom{ for some bounded, ${\cal B}(\Gamma)$-measurable
function $\rho:\Gamma\to\R_+$}}\}\notag
\\&\qquad =\{\mu\}.\label{hjdvhj}\end{align} Let  $G\in D({\cal
E})$ be such that  ${\cal E}(G)=0$. It suffices to prove that
$G=\operatorname{const}$. By the proof of \cite[Lemma~6.1]{AKR4},
without loss of generality, we can suppose that the function $G$
is bounded.

Now, we modify the proof of the part (ii)$\Rightarrow$(iii) of
\cite[Theorem~6.2]{AKR4}. Since $1\in\FC$ and $D^-1=0$, replacing
$G$ by $G-\operatorname{ess\,inf}G$, we may suppose that $G\ge0$,
and, in addition, that $\int G\, d\mu=1$. Define $\nu{:=}G\cdot
\mu$. Since ${\cal E}(G)=0$, by \eqref{hh}, we have that
$G(\gamma\setminus x)-G(\gamma)=0$ $\tilde \mu$-a.e. Since
$\mu\in{\cal G}(z,\phi)$, we have, for any measurable function
$F:\Gamma\times\R^d\to [0,+\infty]$,
\begin{align*}
 \int_\Gamma \nu(d\gamma) \int_{\R^d}\gamma(dx)\, F(\gamma,x)& =\int_\Gamma
\mu(d\gamma)\int_{\R^d}\gamma(dx)\, G(\gamma) F(\gamma,x)\\ &=
\int_\Gamma\mu(d\gamma)\int_{\R^d}\gamma(dx)\, G(\gamma\setminus
x)F(\gamma,x)\\&= \int_\Gamma \mu(d\gamma)\int_{\R^d}z\, dx\,
\exp[-E(x,\gamma)]G(\gamma)F(\gamma\cup x,x)\\&= \int_\Gamma
\nu(d\gamma)\int_{\R^d}z\,dx\, \exp[-E(x,\gamma)]F(\gamma\cup
x,x).
\end{align*} Hence, $\nu\in{\cal G}(z,\phi)$ and, by
\eqref{hjdvhj}, $G=1$ $\mu$-a.e.\quad $\square$

Let us suppose that the potential $\phi$ satisfies the following
condition:

\begin{description}

\item[(LAHT)] ({\it Low activity-high temperature regime})
$$\delta =\int_{\R^d}(1-\exp[-\phi(x)])\, z\,dx <\exp(-1).$$
\end{description}
Under (P) and (LAHT), there exists a unique Gibbs measure
$\mu\in{\cal G}(z,\phi)$, see \cite{Ru69} and
\cite[Theorem~6.2]{Kuna2} (notice that \eqref{fdrtsdrt} and (P)
imply that the correlation functions of any $\mu\in{\cal
G}(z,\phi)$ satisfy the Ruelle bound \eqref{swaswea957} with
$\xi=1$, and hence we can take the constant $C_R$ in
\cite[Theorem~6.2]{Kuna2} to be equal to $e$). The following
statement now  immediately follows from Corollary~\ref{jndfjdj4}.

\begin{cor}[Poincar\'e inequality in the LAHT regime]
Assume that \rom{(P)} and \rom{(LAHT)} are satisfied and consider
the  unique Gibbs measure $\mu\in{\cal G}(z,\phi)$\rom. Then\rom,
\eqref{chdhbbh} holds\rom.
\end{cor}

\begin{center}
{\bf Acknowledgements}\end{center}

We are  grateful to S.~Albeverio, T.~Pasurek, and M.~R\"ockner for
useful discussions. Yu.~K.  gratefully acknowledges the financial
support of the DFG through Projects 436 UKR 113/61 and 436 UKR
113/67. E.~L. gratefully acknowledges the  financial support of
SFB 611, Bonn University.

\end{document}